\documentclass{article}

\usepackage{amssymb,amsmath}
\usepackage[matrix,arrow]{xy}
\usepackage{theorem,graphicx}

\author{Gautam Gangopadhyay}
\title{Counting triangles formula for the first Chern class of a circle bundle}
\date{}

\newtheorem{theorem}{Theorem}

{\theorembodyfont{\rmfamily}    

\newtheorem{definition}{Definition}

}

\def\d{\delta}             \def\:{\colon\,}

\def\sign{\mathop{\rm sign}}

\def\Curv{{\rm Curv}}
\def\Ind{{\rm Ind}}
\let\D\Delta
\def\wD{\widetilde\Delta}

\begin{document}
\maketitle
\begin{abstract}
We consider the problem of the combinatorial computation of the first Chern class of a $S^1$-bundle. N.Mnev found~\cite{Mn} such a formula in terms of canonical shellings. It represents certain invariant of a triangulation computed by analyzing cyclic word in 3-character alphabet associated to the bundle. This curvature is a kind of discretization of Konstevich's curvature differential $2$-form.

We find a new expression of Mnev's curvature by counting triangles in a cyclic word.
This cyclic word represents the multivalued section given by the faces of the triangulation of the fibration.
Our formula is different from that of Mnev. In particular, it is cyclically invariant by its very form. We present also some sample computations of this invariant and also provide a small \emph{Mathematica} code for the computation of this invariant and also We present an independent direct proof of Theorem~\ref{ThGG}. Our arguments are based on the approach developed in~\cite{MK}. The idea of this approach is to consider the simplicial partition of the total space of the fibration as a singular multivalued section that can be used to construct a univalued section. 
\end{abstract}

\section{Introduction}

An \emph{$S^1$-bundle}, or a \emph{circle bundle}, is a locally trivial fiber bundle $\pi:E\to M$ whose fibers are oriented circles. The first Chern class $c_1(\pi)\in H^2(M)$ of an $S^1$-bundle $\pi:E\to M$ is the cohomology class serving as an obstruction to existence of a global section. Assume that the base~$M$ is triangulated and we have a section~$s$ defined over the union of $0$ and $1$-dimensional simplices. Consider a $2$-simplex $\D$. Over this simplex the bundle is trivial, $\pi^{-1}(\D)\approx \D\times S^1$. The section~$s$ restricted to the boundary of $\D$ provides a continuous mapping $s_\D\colon S^1=\partial\D\to S^1$. Define the \emph{index} of that simplex as the degree of the mapping~$s_\D$, that is, the number of rotations of the section~$s$ in the fiber of $\pi$ while the point of the base goes along the boundary circle of~$\D$ (the orientation of~$\partial\D$ is induced from the chosen orientation of~$\D$). The index considered as a function on the set of $2$-simplices can be treated as a simplicial $2$-cochain. This cochain is, in fact, a cocycle, that is, it is closed. By definition, $c_1(\pi)$ is the cohomology class represented by this cocycle.

The defined above cocycle depends on a choice of the section $s$ over the $1$-skeleton of the base.
N.Mnev defined in his paper~\cite{Mn} a different simplicial cocycle for the same cohomology class. This cocycle is independent of any choice and is canonically associated with the combinatorics of a simplicial realization of a given circle bundle. It is worth to mention that Mnev's cocycle is defined over rationals rather than integers.

The combinatorial structure of an $S^1$-bundle is characterized by the so called shelling which is a certain family of cyclic words associated with the simplices of the base. Mnev's formula uses a linear order of letters forming the cyclic word, and the fact that the result is invariant under the cyclic permutation of letters in the word is not obvious. In this paper we derive a different formula for the Mnev's cocycle. Our formula has the advantage that it is cyclically invariant by its very nature. The proof is the direct verification of the equivalence of the two expressions for the corresponding cocycle.

The author is grateful to Professor Maxim Kazarian for the introduction to the subject and valuable discussions.

\section{Combinatorics of a circle bundle}

The topology of a circle bundle can be represented combinatorially as certain additional structure on the simplicial partition of the base. Assume that both the total space and the base are simplicial spaces and the bundle map $\pi$ is a simplicial map. Let $\D$ be a $n$-dimensional simplex of the base. Since the fibers of $\pi$ are one-dimensional, any simplex $\wD$ in the total space that is mapped onto~$\D$ is either $n$-dimensional, or $(n+1)$-dimensional. In the latter case two of the vertices of $\wD$ have a common image. The set of $(n+1)$-simplices in the preimage of~$\D$ is cyclically ordered according to the order of segments of their intersections with a fiber over any internal point of $\D$ in the direction of positive orientation of that fiber. To each of these $(n+1)$-simplices~$\wD$ we can assign a vertex of $\D$, namely, the one whose preimage has two vertices in~$\wD$. Thus we obtain a cyclic word in the alphabet labelled by vertices in $\D$.

The cyclic words associated with different simplices of the base are compatible with one another in the following sense. If $\D'$ is a face of $\D$ (of any dimension), then the cyclic word associated with $\D'$ can be obtained from that for~$\D$ by crossing out all letters corresponding to the vertices in the complement $\D\setminus\D'$. Remark also that each letter appears in the corresponding word at least twice

\begin{definition}
A correspondence that associates to each simplex of a simplicial set a cyclic word in the alphabet labelled by the vertices of the simplex is called a \emph{shelling} if it satisfies the formulated above compatibility condition.
\end{definition}

Thus, to each simplicial realization of an $S^1$-bundle we canonically associate a shelling on the simplicial partition of the base. The correspondence between simplicial realizations of $S^1$-bundles and shellings is bijective: a shelling on a simplicial set defines uniquely both the $S^1$-bundle and the simplicial partition of its total space.

\section{A combinatorial formula due to Mnev}\label{sec3}

In this section we review the formula of N.Mnev~\cite{Mn} for the cocycle representing the first Chern class. The invariance implies that the value of this cocycle on a given $2$-simplex of the base is determined uniquely up to a sign by the cyclic word associated with this simplex. The sign is fixed by the choice of the orientation of the simplex, that is, by a choice of the cyclic order on the set of three letters used in the word.

For a $2$-letter alphabet $a<b$ and a \emph{linear} word $w=w_1w_2\dots w_m$ denote by $k_0$ and $k_1$ the numbers of entrances of the letter~$a$ and~$b$ in the word~$w$, respectively, $k_0+k_1=m$, and by $A\bigsqcup B=\{1,\dots,m\}$ the corresponding splitting of the set of indices: $i\in A\Leftrightarrow w_i=a$, $i\in B\Leftrightarrow w_i=b$. We set
$$\Ind(w)=\frac{\#\{(i,j)\in A\times B\Bigm| i>j\}-\#\{(i,j)\in A\times B\Bigm| i<j\}}{2 k_0 k_1}.$$

Now, consider an alphabet in~$n+1$ letters labelled by the indices $0,1,\dots,n$. For a given word $w=w_1w_2\dots w_m$ in this alphabet we denote by $\d_iw$ the word obtained from~$w$ by crossing out all entrances of the $i$th character of the alphabet. For example, for the word $w=bcabbccacb$ in the alphabet $a<b<c$ we have
\begin{align*}
\d_0(w)&=bcbbcccb &\quad\hbox{in the alphabet }&b<c,\\
\d_1(w)&=caccac   &\quad\hbox{in the alphabet }&a<c,\\
\d_2(w)&=babbab   &\quad\hbox{in the alphabet }&a<b.
\end{align*}

Finally, for a given word in the alphabet in three letters $a<b<c$ denote by $w=w_1w_2\dots w_m$ any of its linear representatives and set
\begin{equation}\label{eq6}
\Curv(w)=\Ind(\d_0w)-\Ind(\d_1w)+\Ind(\d_2w).
\end{equation}

Thus, for the above example $w=bcabbccacb$ we have
\begin{align*}
\Ind(\d_0bcabbccacb)&=\Ind(bcbbcccb)=-\frac{1}{8},\\
\qquad\Ind(\d_1bcabbccacb)&=\Ind(caccac)=0,\\
\qquad\Ind(\d_2bcabbccacb)&=\Ind(babbab)=0, \\
\Curv(bcabbccacb)&=-\frac{1}{8}-0+0=-\frac{1}{8}.
\end{align*}

\begin{theorem}[\cite{Mn}]
The value $\Curv(w)$ is invariant under cyclic permutations of the word~$w$, so it represents a simplicial cochain on the base of a simplicial realization of a circle bundle. This cochain is closed and and represents the cohomology class equal to the first Chern class of the bundle.
\end{theorem}

This theorem is motivated by Kontsevich' construction of piecewise linear connection in a bundle with polygonal fibers~\cite{Kon} but the actual proof given in~\cite{Mn} is not direct: one checks first that $\Curv$ is indeed a cyclically invariant cicycle so that its cohomology class is well defined. Then, some general arguments imply that can not be anything else but the first Chern class up to a factor. The factor is fixed by considering an arbitrary nontrivial example.

\section{Counting triangles formula}
In this section we suggest a new combinatorial formula for the invariant $\Curv$ introduced above. The advantage of our approach is that the suggested formula cyclically invariant by its very form. On the contrary, the terms of the original Mnev's formula~\eqref{eq6} are not cyclically invariant, and the fact that the sum is preserved under cyclic permutations of the word requires a separate proof.

Let $w$ be a cyclic word in a cyclically ordered $3$-character alphabet $(a,b,c)$. Denote by $k_0$, $k_1$, and $k_2$ the numbers of appearances of the characters $a$,~$b$, and~$c$, respectively, in the word~$w$, so that the total length of~$w$ is
$$m = k_0+k_1+k_2.$$

By a \emph{triangle} we mean a choice of three characters in~$w$ such that one of them is~$a$, one is~$b$ and one is~$c$. The cyclic order of letters in a chosen triangle as they appear in the word can agree or disagree with the cyclic order in the alphabet. According to this we call the triangles \emph{oriented} or \emph{disoriented}, respectively. Let~$t^+$ be the number of oriented triangles and~$t^-$ be the number of disoriented ones. So, the total number of possible triangles is
$$ t^+ + t^- = k_0 \times k_1 \times k_2$$

\begin{theorem}\label{ThGG} The invariant $\Curv$ of Section~\ref{sec3} is invariant under cyclic permutations of characters in a word and given explicitly by
$$\Curv(w) = -\frac{t^+ - t^-}{2(t^+ + t^-)}.$$
\end{theorem}

\emph{Proof}. We are going to show that the formula of Theorem~\ref{ThGG} is equivalent to the Mnev's formula~\eqref{eq6}. Let us represent the cyclic word~$w$ as a string (the string is defined up to a cyclic permutation, we just fix any choice): $w=w_1w_2\dots w_m$. Consider the splitting of set of indices $\{1,\dots,m\}=A\bigsqcup B\bigsqcup C$ so that $w_i=a\Leftrightarrow i\in A$, and similarly for $B$ and $C$. We have: $k_0=|A|$, $k_1=|B|$, $k_2=|C|$.

By the Mnev's formula,
\begin{align*}
\Curv(w) &= \Ind(\d _{0}w)-\Ind(\d _{1}w) + \Ind(\d _{2}w)\\
         &= -\frac{\!\!\sum\limits_{j\in B,k \in C}\!\!\sign(k-j)}{2\,k_1\, k_2}
            +\frac{\!\!\sum\limits_{i\in A,k \in C}\!\!\sign(k-i)} {2\,k_0\, k_2}
            -\frac{\!\!\sum\limits_{i\in A,j \in B}\!\!\sign(j-i)} {2\,k_0\, k_1}\\
         &= -\frac{\!\!\!\sum\limits_{i\in A,j \in B,k \in C}\!\!\!S(i,j,k)} {2\,k_0\, k_1 \, k_2},
\end{align*}
where
$$S(i,j,k)=\sign(k-j)-\sign(k-i)+\sign(j-i).$$
The value of $S(i,j,k)$ depends on the linear order of the indices $i$,~$j$, and~$k$. There are 6 possible such orders, and the computation of the value $S(i,j,k)$ is given in the following table
$$\begin{array}{|c|c|}\hline
   \text{ordering}&S(i,j,k)\\\hline
 i<j<k  &  (+1)-(+1)+(+1) = (+1)\\
 j<k<i  &  (-1)-(-1)+(+1) = (+1)\\
 k<i<j  &  (-1)-(-1)+(+1) = (+1)\\
 i<k<j  &  (-1)-(+1)+(+1) = (-1)\\
 j<i<k  &  (+1)-(+1)+(-1) = (-1)\\
 k<j<i  &  (-1)-(-1)+(-1) = (-1)\\\hline
 \end{array}$$

 One can see that the value of $S(i,j,k)$ is $+1$ if the cyclic order of the indices $i,j,k$ is positive and the value is $-1$ otherwise. It remains to observe that the triples $(i,j,k)$, $i\in A$, $j\in B$, $k\in C$ enumerate triangles, and the number of positive and negative values of $S(i,j,k)$ are exactly~$t^+$ and~$t^-$, respectively. This proves the formula of Theorem~\ref{ThGG}.

\bigskip
\bigskip
A short \emph{Mathematica} code implementing the computation of $\Curv$ from  the formula of Theorem~\ref{ThGG} is presented below.

\bigskip
\noindent
\includegraphics[scale=0.4]{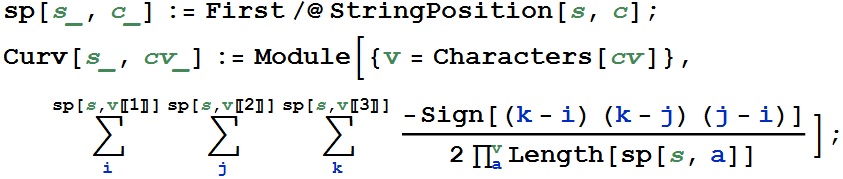}\\ \\
\verb|Curv["bcabbccacb", "abc"]|
\\\verb|-1/8|
\\
\\\verb|Curv["ddbbccdbc", "cbd"]|
\\\verb|5/18|
\\
\\\verb|Curv["papaspaspsa", "aps"]|
\\\verb|1/24|

\section{Independent Proof}

In this section, we present an independent direct proof of Theorem~\ref{ThGG}. Our arguments are based on the approach developed in~\cite{MK}. The idea of this approach is to consider the simplicial partition of the total space of the fibration as a singular multivalued section that can be used to construct a univalued section. First, we constract a section over $0$-simplices, then, over $1$-simplices, then, the extension to $2$-simplices meets an obstruction which is a cochain representing the first Chern class, as in Introduction. There is a freedom in a choice of a section, so we consider all possible choices, and then take the average. This eliminates the ambiguity in a choice of a section and leads to the statement of Theorem~\ref{ThGG}.

Now we explain this procedure in more details. At the vertices, we choose the value of the section \emph{at the middle point of one of the segments forming the simplicial partition of the fiber}. This section has a natural extension to all neighboring simplices by setting its value also at the middle of the segments into which simplices of the total space split the fiber of the bundle.

Then, over the $1$-simplices we have two chosen sections, one coming from one end of the $1$-simplex and one from another end. These sections do not intersect each other, and we have to connect them to form a single continues section over the segment by some arbitrary way. We do this along one of the two arcs of the circle either in the clockwise ($C_S$) or in the counterclockwise ($AC_S$) direction. A choice of one of these two arcs is another freedom that we have to take into account.

Consider a simplex of the base with vertices $a,b,c$. Assume that the chosen orientation of this simplex corresponds to the cyclic order $abc$ of its vertices. The combinatorial structure of the fibration over this simplex corresponds to a cyclic word with $k_0,k_1,k_2$ appearances of the letters $a,b,c$. We have got, therefore, $8k_0k_1k_2$ choices for a section along its boundary: there is $k_0$ choices to choose a section over $a$, etc, and there are two possibilities to extend this section along the segments of the triangle.

First, we fix a choice of the values of the sections over the vertices, which is equivalent to a choice of one of $k_0k_1k_2$ triangles, and then, we take the average over $8$ possible choices of the arcs. For a fixed triangle the average gives $\pm1/2$ as shown below, the sign depends on the orientation of the triangle. This leads to the required formula.

Let $A$, $B$, and $C$ be the chosen values of the section over the vertices $a$, $b$, and $c$, respectively. Assume that these points are situated on the fiber in the cyclic order $ACB$. We compute the index of each of $8$ extensions of the section to the sides of the triangle $abc$ by counting the rotation number of the obtained path connecting the points $A$, $B$, and $C$. The computation of these indices is presented in the table below. In this table, we denote by $C_S$ the choice of the arc in the clockwise direction and by $AC_S$ the choice in the counterclockwise direction
\def\ar#1{{}\stackrel{#1}{\longrightarrow}{}}
$$
\begin{array}{@{A}c@{B}c@{C}c@{A\quad}|c}
\multicolumn{3}{c|}{\hbox{Section}}&\hbox{Index}\\\hline
\ar{C_S}&\ar{C_S}&\ar{C_S}&-1\\
\ar{C_S}&\ar{C_S}&\ar{AC_S}&0\\
\ar{C_S}&\ar{AC_S}&\ar{AC_S}&+1\\
\ar{C_S}&\ar{AC_S}&\ar{C_S}&0\\
\ar{AC_S}&\ar{AC_S}&\ar{AC_S}&+2\\
\ar{AC_S}&\ar{AC_S}&\ar{C_S}&+1\\
\ar{AC_S}&\ar{C_S}&\ar{AC_S}&+1\\
\ar{AC_S}&\ar{C_S}&\ar{C_S}&0
\end{array}
$$

Averaging these values we get
$$\frac{-1+0+1+0+2+1+1+0}{8}=+\frac{1}{2}.$$

Similar computations show that in the case when the points $A$, $B$, and $C$ go the cyclic order $ABC$ on the fiber the average index is $-1/2$.

Therefore, the total average index of the possible choices of the section along the boundary of the simplex is
\begin{multline*}
\frac{1}{k_0k_1k_2}\sum_{\hbox{Choices of triangles}}\frac{1}{8}\sum_{\hbox{Choices of arcs}}\hbox{index of a chosen section}\\
=\frac{1}{k_0k_1k_2}\Bigl(\frac{1}{2}t^- - \frac{1}{2}t^+\Bigr)=\frac{t^--t^+}{2\,k_0k_1k_2}.
\end{multline*}
This computation completes the proof.

\section{Conclusion}

We proved our formula of curvature by counting triangles in a cyclic word is equivalent to the Mnev' formula of the problem of the combinatorial computation of the first Chern class of a $S^1$-bundle. This formula suggests its generalizations to higher dimensions. There is a possible formula for the powers of the first Chern Class is discussed in~\cite{MNEV}. we presented an independent direct proof of Theorem~\ref{ThGG} by using the technique of multivalued sections introduced in~\cite{MK}.

\end{document}